\documentclass{amsart}
\usepackage{amsfonts}
\usepackage{amssymb}
\usepackage{graphicx}
\usepackage{xypic}
\usepackage{psfrag}

\newcommand{\N}{{\mathbb N}}
\newcommand{\Z}{{\mathbb Z}}

\newcommand{\Th}{{\text{th}}}
\newcommand{\St}{{\text{st}}}
\newcommand{\partof}{\,\vdash\,}
\newcommand{\Schub}{{\mathfrak{S}}}
\newcommand{\bull}{{\sssize \bullet}}
\newcommand{\rank}{\operatorname{rank}}

\newcommand{\pic}[2]{\includegraphics[scale=0.#1]{#2.eps}}

\newcommand{\picC}[1]{\includegraphics[scale=0.50]{#1.eps}}

\newcommand{\smpicskip}{\vspace{0.1cm}}

\theoremstyle{plain}

\newtheorem*{thm}{Theorem}

\newtheorem*{cor}{Corollary}

\theoremstyle{definition}
\newtheorem{example}{Example}

\newcommand{\refeqn}[1]{(\ref{#1})}

\newcommand{\lab}[1]{\label{#1}}

\newcommand{\partit}[1]{\raisebox{-5pt}{\includegraphics[scale=.60]{#1.eps}}}

\newcommand{\Pa}{\partit{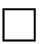}}

\newcommand{\Pb}{\partit{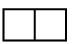}}
\newcommand{\Paa}{\partit{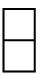}}

\newcommand{\Pba}{\partit{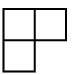}}
\newcommand{\Paaa}{\partit{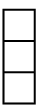}}

\newcommand{\Pca}{\partit{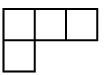}}

\newcommand{\Pbaa}{\partit{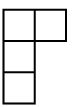}}

\psfrag{rij-1-rij}{$r_{i,j-1} - r_{ij}$}
\psfrag{ri+1j-rij}{$r_{i+1,j} - r_{ij}$}
\psfrag{R01}{$R_{01}$}
\psfrag{R12}{$R_{12}$}
\psfrag{Rn-1n}{$R_{n-1,n}$}
\psfrag{sig1}{$\sigma_1$}
\psfrag{sig2}{$\sigma_2$}
\psfrag{tau1}{$\tau_1$}
\psfrag{taun-1}{$\tau_{n-1}$}
\psfrag{Ri-1i}{$R_{i-1,i}$}
\psfrag{sigi}{$\sigma_i$}
\psfrag{taui-1}{$\tau_{i-1}$}
\psfrag{Rectangle diagram for r}{Rectangle diagram for $r$}
\psfrag{Rectangle diagram for r'}{Rectangle diagram for $r'$}
\psfrag{Rectangle diagram for r''}{Rectangle diagram for $r''$}
\psfrag{Empty}{empty}
\psfrag{empty}{empty}
\psfrag{rectangles}{rectangles}
\psfrag{Rectangles}{rectangles}
\psfrag{for u}{for $u$}
\psfrag{for w}{for $w$}

\title{Stanley symmetric functions and quiver varieties}
\date{\today}
\author{Anders Skovsted Buch}
\address{Massachusetts Institute of Technology \\
  Building 2, Room 248 \\
  77 Massachusetts Avenue \\
  Cambridge, MA 02139
}
\email{abuch@math.mit.edu} 

\begin{document}
\maketitle

\section{Introduction}

The purpose of this paper is to show a connection between Stanley
symmetric functions and the formula for quiver varieties given in
\cite{buch.fulton:chern}.  Recall that a simple reflection in the
symmetric group $S_m$ is a transposition that interchanges two
consecutive integers.  A {\em reduced word\/} for a permutation $w \in
S_m$ is a tuple of simple reflections $(\tau_1, \tau_2, \dots,
\tau_\ell)$ with $\ell = \ell(w)$ the length of $w$, such that $w =
\tau_1 \tau_2 \cdots \tau_\ell$.  Stanley asked how many reduced words
does a permutation $w$ have.

To answer this question, Stanley \cite{stanley:on*1} defined a power
series $F_w(x)$ in infinitely many variables $x_1, x_2, \dots$; it is
homogeneous of degree $\ell = \ell(w)$, has non-negative integer
coefficients, and the number of reduced words for $w$ is the
coefficient in $F_w(x)$ of the monomial $x_1 x_2 \cdots x_\ell$.
Stanley then proved that this power series is symmetric.  This implies
that it can be written in the basis of Schur functions:
\begin{equation}
\lab{eqn:stanley}
F_w(x) = \sum_{\lambda \partof \ell} 
   \alpha_{w\lambda} \, s_\lambda(x) 
\end{equation}
where the sum is over all partitions $\lambda$ of $\ell$ and the
coefficients $\alpha_{w \lambda}$ are integers.  Since the coefficient
of $x_1 x_2 \cdots x_\ell$ in a Schur function $s_\lambda(x)$ is equal
to the number $f^\lambda$ of standard Young tableaux of shape
$\lambda$ (see e.g.\ \cite{macdonald:symmetric*1} or
\cite{fulton:young}), it follows that the number of reduced words for
$w$ is given as
\[ \sum_{\lambda \partof \ell} \alpha_{w \lambda} \,f^\lambda \,. \]
The constants $f^\lambda$ are considered well understood, so only a
description of the coefficients $\alpha_{w \lambda}$ remained to be
found.  Stanley credits Edelman and Greene for proving that these
coefficients are non-negative \cite{edelman.greene:balanced} (see also
\cite{lascoux.schutzenberger:structure}).  Fomin and Greene have shown
that $\alpha_{w \lambda}$ is equal to the number of semistandard Young
tableaux $T$ of shape $\lambda$, such that the column word of $T$ is a
reduced word for $w$ \cite{fomin.greene:noncommutative}.  Another
useful fact is that $\alpha_{w^{-1}\lambda} = \alpha_{w\lambda'}$
where $\lambda'$ is the conjugate of $\lambda$
\cite{lascoux.schutzenberger:structure},
\cite[(7.22)]{macdonald:notes}.

Stanley's symmetric function is known to be a limit of Schubert
polynomials $\Schub_w(x)$ defined by Lascoux and Sch\"utzenberger
\cite{lascoux.schutzenberger:polynomes}, \cite{macdonald:notes}.  For
$n \in \N$, let $1^n \times w \in S_{n+m}$ denote the shifted
permutation which acts as the identity on $1,\dots,n$ and maps $i$ to
$w(i-n)+n$ for $n+1 \leq i \leq n+m$.  If one specializes to finitely
many variables $x_1,x_2, \dots, x_N$, then
\begin{equation}
\lab{eqn:stanstab}
 F_w(x_1,\dots,x_N,0,0,\ldots) = 
   \Schub_{1^n \times w^{-1}}(x_1,\dots,x_N,0,0,\ldots)
\end{equation}
for all sufficiently large $n$ \cite[(7.18)]{macdonald:notes}.

The formula for quiver varieties given in \cite{buch.fulton:chern}
specializes to a formula for the double Schubert polynomial
$\Schub_w(x;y)$ for the permutation $w \in S_m$:
\begin{equation}
\lab{eqn:schub}
\Schub_w(x;y) = \sum c_w(a,b,\lambda) \, 
   y_2^{a_2} \cdots y_{m-1}^{a_{m-1}} \, 
   (-x_2)^{b_2} \cdots (-x_{m-1})^{b_{m-1}} \,
   s_\lambda(x/y) \,.
\end{equation}
The sum is over exponents $a_2,\dots,a_{m-1}$ and $b_2,\dots,b_{m-1}$
and a partition $\lambda$.  The coefficients $c_w(a,b,\lambda)$ are
special cases of a large class of generalized Littlewood-Richardson
coefficients, which are conjectured to be non-negative and given by a
generalized Littlewood-Richardson rule \cite{buch.fulton:chern}.

The main result in this paper is that the coefficient
$c_w(0,0,\lambda)$ (corresponding to zero exponents) is equal to
Stanley's coefficient $\alpha_{w^{-1} \lambda}$.  In this way,
(\ref{eqn:schub}) writes a Schubert polynomial $\Schub_w(x)$ as a
symmetric polynomial equal to Stanley's symmetric function for
$w^{-1}$ plus a non-symmetric polynomial.

In Section \ref{sec:formula} and Section \ref{sec:schubert} we review
the results of \cite{buch.fulton:chern}.  In Section \ref{sec:stable}
we prove the identity $\alpha_{w^{-1} \lambda} = c_w(0,0,\lambda)$ and
use this to give a new proof of Stanley's result \cite{stanley:on*1}
that the symmetric function $F_{w_0}(x)$ for the longest permutation
$w_0$ in $S_m$ is equal to the Schur function $s_\lambda(x)$ for the
staircase partition $\lambda = (m-1,m-2,\dots,1)$.  In Section
\ref{sec:redundant} we use geometry of degeneracy loci to prove a
generalization of the well known formula for a Schubert polynomial of
a product of two permutations.  Finally, in Section \ref{sec:conj}, we
discuss the relations to the conjectured generalized
Littlewood-Richardson rule in
\cite{buch.fulton:chern}.

We thank S.~Fomin, W.~Fulton, and F.~Sottile for helpful discussions.

\section{A formula for quiver varieties}
\lab{sec:formula}

Let $X$ be a non-singular complex variety and $E_0 \to E_1 \to E_2 \to
\cdots \to E_n$ a sequence of vector bundles and vector bundle maps
over $X$.  A set of {\em rank conditions\/} for this sequence is a
collection $r = (r_{ij})$ of non-negative integers, for $0 \leq i < j
\leq n$.  Let $\Omega_r(E_\bull) \subset X$ be the locus where each
composite map $E_i \to E_j$ has rank at most $r_{ij}$:
\[ \Omega_r(E_\bull) = \{ x \in X \mid \rank(E_i(x) \to E_j(x)) \leq
   r_{ij} ~\forall i < j \} \,.
\]
This locus determines a cohomology class $[\Omega_r(E_\bull)]$ in the
cohomology ring $H^*(X) = H^*(X;\Z)$ of $X$.

For convenience we set $r_{ii} = \rank(E_i)$.  We will require that
the rank conditions $r = (r_{ij})$ can {\em occur\/}, i.e.\ there should exist
a sequence of vector spaces and linear maps $V_0 \to V_1 \to \cdots
\to V_n$ such that $\dim(V_i) = r_{ii}$ and $\rank(V_i \to V_j) =
r_{ij}$ for all $i < j$.  This is equivalent to demanding $r_{ij} \leq
\min(r_{i,j-1}, r_{i+1,j})$ for $i < j$ and 
$r_{ij} - r_{i,j-1} - r_{i+1,j} + r_{i+1,j-1} \geq 0$
for $j-i \geq 2$ \cite{abeasis.del-fra:degenerations}.

If the locus $\Omega_r(E_\bull)$ is not empty, its codimension is at
most 
\[ d(r) = \sum_{i < j} (r_{i,j-1} - r_{ij})(r_{i+1,j} - r_{ij}) \,. \]
Furthermore, this codimension is obtained for generic choices of
bundle maps $E_i \to E_{i+1}$ \cite{buch.fulton:chern},
\cite{lakshmibai.magyar:degeneracy},
\cite{abeasis.del-fra.ea:geometry}.  The main result of
\cite{buch.fulton:chern} gives a formula for the cohomology class
of the degeneracy locus $\Omega_r(E_\bull)$, when it has this expected
codimension $d(r)$.

To explain this formula, we need some notation.  Let $\Lambda =
\Z[h_1, h_2, \ldots]$ be the ring of symmetric functions.  The
variable $h_i$ can be identified with the complete symmetric function
of degree $i$ in variables $x_i$, $i \in \N$.  If $I = (a_1, a_2,
\dots, a_p)$ is a sequence of integers, define the Schur function $s_I
\in \Lambda$ to be the determinant of the $p \times p$ matrix whose
$(i,j)^\Th$ entry is $h_{a_i+j-i}$:
\[ s_I = \det(h_{a_i+j-i})_{1\leq i,j \leq p} \,. \]
(Here one sets $h_0 = 1$ and $h_{-q} = 0$ for $q > 0$.)  A Schur
function is always equal to either zero, or plus or minus a Schur
function for a partition $\lambda$:
\[ s_I = \begin{cases} 0 \\ \pm s_\lambda \,. \end{cases} \]
This follows from interchanging the rows of the matrix defining $s_I$.
Furthermore, the Schur functions given by partitions form a basis for
the ring of symmetric functions \cite{macdonald:symmetric*1},
\cite{fulton:young}.

If $E$ is a vector bundle over $X$, let $c_i(E) \in H^{2i}(X)$ denote
its $i^\Th$ Chern class.  Given two vector bundles $E$ and $F$ of
ranks $e$ and $f$ over $X$ one can define a ring homomorphism
\[ \Lambda \to H^*(X) \]
which maps $h_i$ to the coefficient of $t^i$ in the formal power
series expansion of the quotient
\[ \frac{c_t(E^\vee)}{c_t(F^\vee)} = 
   \frac{1 - c_1(E) t + \dots + (-1)^e c_e(E) t^e}
   {1 - c_1(F) t + \dots + (-1)^f c_f(F) t^f} \,. 
\]
We let $s_\lambda(F - E) \in H^*(X)$ denote the image of
$s_\lambda$ by this map.  (If $E$ and $F$ have Chern roots $x_1,
\dots, x_e$ and $y_1, \dots, y_f$ respectively, the notation
$s_\lambda(y/x) = s_\lambda(F - E)$ is also common.)

When the degeneracy locus $\Omega_r(E_\bull)$ has its expected
codimension $d(r)$, its cohomology class is equal to a linear
combination of products of Schur polynomials in differences of
consecutive bundles:
\[ [\Omega_r(E_\bull)] = \sum_\mu c_\mu(r) \, s_{\mu_1}(E_1 - E_0) 
   s_{\mu_2}(E_2 - E_1) \cdots s_{\mu_n}(E_n - E_{n-1}) \,. 
\]
Here the sum is over sequences of partitions $\mu = (\mu_1, \mu_2, \dots,
\mu_n)$.  The coefficients $c_\mu(r)$ are integer constants depending
on the rank conditions and the sequence $\mu$.  They are determined by 
a combinatorial algorithm which we will describe next.

Start by arranging the rank conditions $r = (r_{ij})$ in a {\em rank
diagram}:
\[ \begin{matrix}
E_0 & \to & E_1 & \to & E_2 & \to & \cdots & \to & E_n 
\vspace{0.1cm} \\
r_{00} && r_{11} && r_{22} && \cdots && r_{nn} \\
& r_{01} && r_{12} && \cdots && r_{n-1,n} \\
&& r_{02} && \cdots && r_{n-2,n} \\
&&& \ddots \\
&&&& r_{0n}
\end{matrix} \]
In this diagram, replace each small triangle of numbers
\[ \begin{matrix}
r_{i,j-1} && r_{i+1,j} \\
& r_{ij}
\end{matrix} \]
by a rectangle $R_{ij}$ with $r_{i+1,j} - r_{ij}$ rows and $r_{i,j-1}
- r_{ij}$ columns.
\[ R_{ij} = \raisebox{-18pt}{\picC{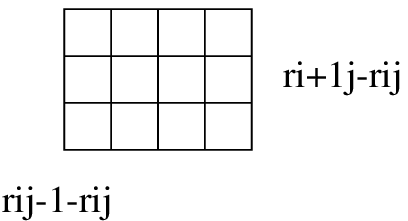}} \]
These rectangles are arranged in a {\em rectangle diagram}:
\[ \begin{matrix}
R_{01} && R_{12} && \cdots && R_{n-1,n} \\
& R_{02} && \cdots && R_{n-2,n} \\
&& \ddots \\
&&& R_{0n}
\end{matrix} \]
The information carried by the rank conditions is very well
represented in this diagram.  First, the expected codimension $d(r)$
for the locus $\Omega_r(E_\bull)$ is equal to the total number of
boxes in the rectangle diagram.  Furthermore, the condition that the
rank conditions can occur is equivalent to saying that the rectangles
get narrower when one travels south-west, while they get shorter when
one travels south-east.  Finally, the algorithm that computes the
coefficients $c_\mu(r)$ depends only on the rectangle diagram.

We will define this algorithm by constructing an element $P_r$ in the
$n^\Th$ tensor power of the ring of symmetric functions
$\Lambda^{\otimes n}$, so that
\[ P_r = 
   \sum_\mu c_\mu(r) \, s_{\mu_1} \otimes \cdots \otimes s_{\mu_n} \,. 
\]
This is done by induction on $n$.  When $n = 1$ (corresponding to a
sequence of two vector bundles), the rectangle diagram has only one
rectangle $R = R_{01}$.  In this case we set
\[ P_r = s_R \in \Lambda^{\otimes 1} \,, \]
where $R$ is identified with the partition for which it is the Young
diagram.  This case recovers the Giambelli-Thom-Porteous formula.

If $n \geq 2$, we let $\Bar r$ denote the bottom $n$ rows of the rank
diagram.  Then $\Bar r$ is a valid set of rank conditions, so by
induction we can assume that 
\begin{equation}
\lab{eqn:prbar}
P_{\Bar r} = \sum_\mu c_\mu(\Bar r) \,
s_{\mu_1} \otimes \cdots \otimes s_{\mu_{n-1}}
\end{equation}
is a well defined element of $\Lambda^{\otimes n-1}$.  Now $P_r$ is
obtained from $P_{\Bar r}$ by replacing each basis element $s_{\mu_1}
\otimes \cdots \otimes s_{\mu_{n-1}}$ in (\ref{eqn:prbar}) with the
sum 
\[ \sum_{\stackrel{\sigma_1,\dots,\sigma_{n-1}}{\tau_1,\dots,\tau_{n-1}}}
   \left(\prod_{i=1}^{n-1} c^{\mu_i}_{\sigma_i \tau_i}\right)
   s_{\picC{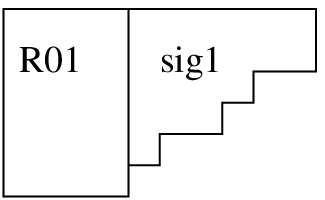}} \otimes \cdots \otimes
   s_{\picC{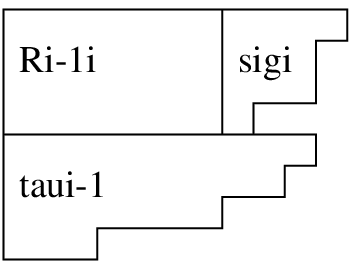}} \otimes \cdots \otimes
   s_{\picC{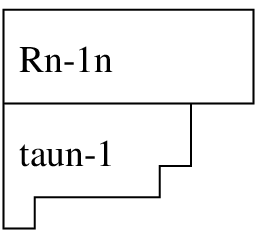}} \,.
\]
This sum is over all partitions $\sigma_1, \dots, \sigma_{n-1}$ and
$\tau_1, \dots, \tau_{n-1}$ such that $\sigma_i$ has fewer rows than
$R_{i-1,i}$ and each Littlewood-Richardson coefficient
$c^{\mu_i}_{\sigma_i \tau_i}$ is non-zero.  A diagram consisting of a
rectangle $R_{i-1,i}$ with (the Young diagram of) a partition
$\sigma_i$ attached to its right side, and $\tau_{i-1}$ attached
beneath should be interpreted as the sequence of integers giving the
number of boxes in each row of this diagram.

It can happen that the rectangle $R_{i-1,i}$ is empty, since the
number of rows or columns can be zero.  If the number of rows is zero,
then $\sigma_i$ is required to be empty, and the diagram is the Young
diagram of $\tau_{i-1}$.  If the number of columns is zero, then the
algorithm requires that the length of $\sigma_i$ is at most equal to
the number of rows $r_{ii} - r_{i-1,i}$, and the diagram consists of
$\sigma_i$ in the top $r_{ii} - r_{i-1,i}$ rows and $\tau_{i-1}$ below
this, possibly with some zero-length rows in between.

\section{Schubert polynomials}
\lab{sec:schubert}

Let $w \in S_{m+1}$ be a permutation, and let $E_\bull$ be a sequence
of bundles over $X$
\[ F_1 \subset F_2 \subset \cdots \subset F_m \to 
    G_m \twoheadrightarrow G_{m-1} \twoheadrightarrow \cdots 
    \twoheadrightarrow G_1
\]
consisting of a full flag with a general map to a dual full flag.
Define the locus
\[ \Omega_w = \{ x \in X \mid \rank(F_q(x) \to G_p(x)) \leq r_w(p,q)
   ~\forall p,q \}
\]
where $r_w(p,q) = \# \{ i \leq p \mid w(i) \leq q \}$.  Fulton has
proved \cite{fulton:flags} that the cohomology class of this locus is
given by the double Schubert polynomial defined by Lascoux and
Sch\"utzenberger \cite{lascoux.schutzenberger:polynomes}:
\[ [\Omega_w] = \Schub_w(x_1,\dots,x_m ; y_1,\dots,y_m) \]
where $x_i = c_1(\ker(G_i \to G_{i-1}))$ and $y_i = c_1(F_i/F_{i-1})$.
Now $\Omega_w = \Omega_r(E_\bull)$ where $r = (r_{ij})$ are the
obvious rank conditions.  This means that the double Schubert
polynomial becomes a special case of the quiver formula:
\[ \begin{split}
\Schub_w(x;y) &= [\Omega_r(E_\bull)] \\
&= \sum c_\mu(r)\, s_{\mu_1}(F_2 - F_1) \cdots
   s_{\mu_{m-1}}(F_m - F_{m-1}) \cdot s_{\mu_m}(G_m - F_m) \cdot \\
& \hspace{2cm} s_{\mu_{m+1}}(G_{m-1} - G_m) \cdots s_{\mu_{2m-1}}(G_1 - G_2)
\end{split} \]
As noted in \cite{buch.fulton:chern}, significant simplifications
can be made by using the equalities
\[ s_\lambda(F_{i+1} - F_i) = s_\lambda(y_{i+1}) = 
   \begin{cases} 
   y_{i+1}^a & \text{if $\lambda = (a)$ is a row with $a$ boxes} \\
   0 & \text{otherwise}
   \end{cases}
\]
and
\[ s_\lambda(G_i - G_{i+1}) = s_\lambda(0/x_{i+1}) =
   \begin{cases}
   (-x_{i+1})^b & \text{if $\lambda = (1^b)$ is a column with $b$
   boxes} \\
   0 & \text{otherwise.}
   \end{cases}
\]
Using this and the fact that $s_\lambda(G_m - F_m)$ is the
super-symmetric Schur polynomial $s_\lambda(x/y)$ in the variables
$x_1,\dots,x_m$ and $y_1,\dots,y_m$, we obtain a formula
\begin{equation}
\lab{eqn:newschub}
\Schub_w(x;y) = \sum c_w(a,b,\lambda) \, y_2^{a_2} \cdots y_m^{a_m}
   \, (-x_2)^{b_2} \cdots (-x_m)^{b_m} \, s_\lambda(x/y) \,.
\end{equation}
The sum is over exponents $a_2,\dots,a_m$ and $b_2,\dots,b_m$, and a
single partition $\lambda$, and $c_w(a,b,\lambda)$ is the coefficient
$c_\mu(r)$ for the sequence of partitions 
\[ \mu = ((a_2), \dots, (a_m), \lambda, 
   (1^{b_m}), \dots, (1^{b_2})) \,. 
\]

\begin{example}
\lab{exm:formula}
For the permutation $w = 2\,4\,3\,1$ we get the rank diagram
\setcounter{MaxMatrixCols}{11}
\[ \begin{matrix}
F_1 & \subset & F_2 & \subset & F_3 & \to & G_3 & \twoheadrightarrow &
G_2 & \twoheadrightarrow & G_1 
\vspace{0.1cm} \\
1 && 2 && 3 && 3 && 2 && 1 \\
& 1 && 2 && 2 && 2 && 1 \\
&& 1 && 1 && 1 && 1 \\
&&& 0 && 1 && 1 \\
&&&& 0 && 1 \\
&&&&& 0
\end{matrix} \]
which in turn gives the rectangle diagram:
\[ \pic{45}{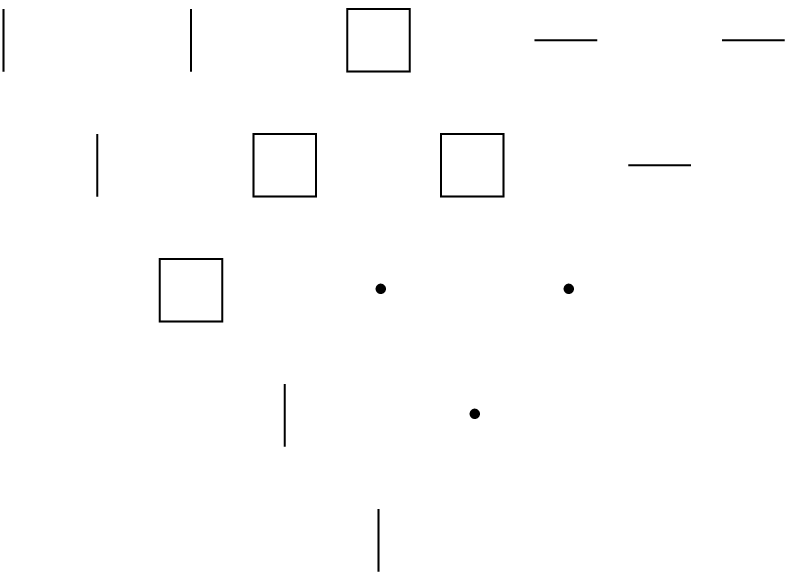} \]
The bottom three rows of this rectangle diagram gives 
\[ P_{\Bar{\Bar r}} = s_{\Pa} \otimes 1 \otimes 1 \,; \]
using the algorithm we then get
\[ P_{\Bar r} = s_{\Pa} \otimes s_{\Pa} \otimes s_{\Pa} \otimes 1 ~+~
   1 \otimes s_{\Paa} \otimes s_{\Pa} \otimes 1
\]
and 
\[ \begin{split}
P_r =~ &
  s_{\Pa} \otimes s_{\Pa} \otimes s_{\Pb} \otimes 1 \otimes 1 ~+~
  s_{\Pa} \otimes s_{\Pa} \otimes s_{\Pa} \otimes s_{\Pa} \otimes 1 ~+ \\
&  s_{\Pa} \otimes 1 \otimes s_{\Pba} \otimes 1 \otimes 1 ~+~
  s_{\Pa} \otimes 1 \otimes s_{\Paa} \otimes s_{\Pa} \otimes 1 ~+ \\
&  1 \otimes s_{\Paa} \otimes s_{\Pb} \otimes 1 \otimes 1 ~+~ 
  1 \otimes s_{\Paa} \otimes s_{\Pa} \otimes s_{\Pa} \otimes 1 ~+ \\
&  1 \otimes s_{\Pa} \otimes s_{\Pba} \otimes 1 \otimes 1 ~+~ 
  1 \otimes s_{\Pa} \otimes s_{\Paa} \otimes s_{\Pa} \otimes 1 ~+ \\
&  1 \otimes 1 \otimes s_{\Pbaa} \otimes 1 \otimes 1 ~+~
  1 \otimes 1 \otimes s_{\Paaa} \otimes s_{\Pa} \otimes 1 \,.
\end{split} \]
This gives the formula
\begin{multline*} 
\Schub_w(x;y) =
  y_2\, y_3\, s_{\Pb}(x/y) - x_3\, y_2\, y_3\, s_{\Pa}(x/y)
  + y_2\, s_{\Pba}(x/y) 
  - x_3\, y_2\, s_{\Paa}(x/y) \\
  + y_3\, s_{\Pba}(x/y) - x_3\, y_3\, s_{\Paa}(x/y) 
  + s_{\Pbaa}(x/y) - x_3\, s_{\Paaa}(x/y) \,.
\end{multline*}
\end{example}

In general, the rectangle diagram associated to a permutation $w \in
S_{m+1}$ contains only empty rectangles and $1 \times 1$ rectangles,
and all of the non-empty ones are located in a diamond below the
rectangle $R_{m-1,m}$.  In other words, if $R_{ij}$ is not empty then
$i \leq m-1$ and $j \geq m$.  In fact, $R_{ij}$ is non-empty if and
only if the diagram $D'(w)$ from \cite{fulton:flags} has a box in
position $(2m-j, i+1)$, and this happens exactly when $w(2m+1-j) \leq
i+1$ and $w^{-1}(i+2) \leq 2m-j$ \cite{buch.fulton:chern}.

\section{Stable Schubert polynomials}
\lab{sec:stable}

In this section we will apply the quiver formula for Schubert
polynomials to calculate Stanley symmetric functions.  Let $w \in
S_{m+1}$ be a permutation and $r = (r_{ij})$ the corresponding rank
conditions.  Notice at first that the rank diagram for the one step
shifted permutation $1 \times w$ is obtained by adding one to each
number $r_{ij}$ in the rank diagram for $w$, and putting an extra row
of ones on the sides of this diagram.  For example, if $w = 3\,1\,2$,
this looks like: \smpicskip
\[
\begin{matrix}
1 && 2 && 2 && 1 \\
& 1 && 1 && 1 \\
& & 1 && 0 \\
& & & 0 \\
\\ \mbox{}
\end{matrix} 
\hspace{.5cm} \rightsquigarrow \hspace{.5cm}
\begin{matrix}
1 && 2 && 3 && 3 && 2 && 1 \\
& 1 && 2 && 2 && 2 && 1 \\
& & 1 && 2 && 1 && 1 \\
& & & 1 && 1 && 1 \\
& & & & 1 && 1 \\
& & & & & 1
\end{matrix}
\]
This means that the rectangle diagram for $1 \times w$ is obtained by
adding a rim of empty rectangles to the sides of the rectangle diagram
for $w$.
\smpicskip
\[ \raisebox{-10pt}{\pic{45}{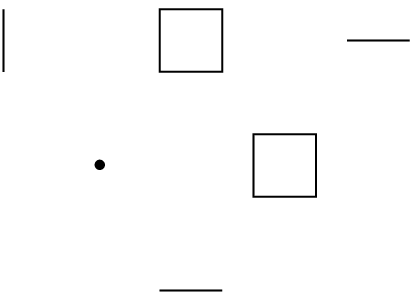}}
   \hspace{.5cm} \rightsquigarrow \hspace{.5cm}
   \raisebox{-43pt}{\pic{45}{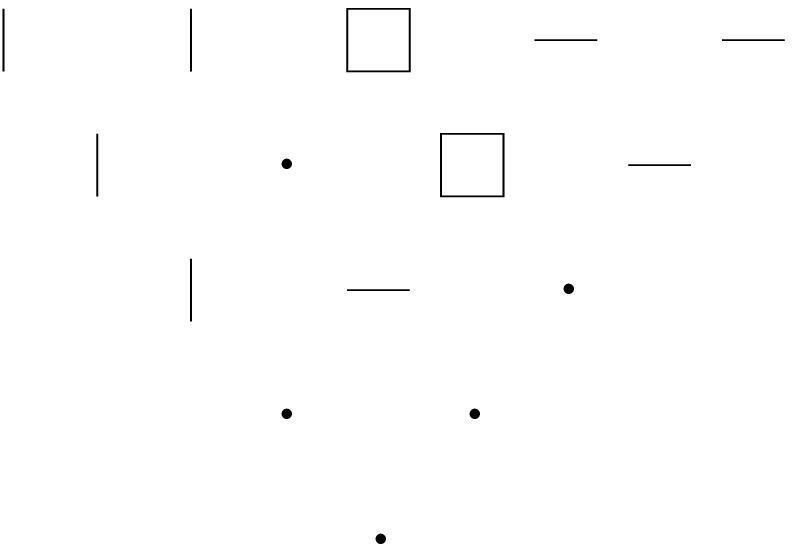}}
   \smpicskip
\]
Similarly, one obtains the rectangle diagram for $1^n \times w$ by
adding $n$ rims of empty rectangles to the rectangle diagram for $w$.

Let $P_r \in \Lambda^{\otimes 2m-1}$ be the element associated to the
rank conditions $r = (r_{ij})$ for $w$.  The above comparison of
rectangle diagrams then shows that $1^n \times w$ corresponds to the
element 
\[ \underbrace{1 \otimes \cdots \otimes 1}_n \otimes P_r \otimes 
   \underbrace{1 \otimes \cdots \otimes 1}_n
   \in \Lambda^{\otimes 2m+2n-1} \,.
\]
By \refeqn{eqn:newschub} this gives us
\[
 \Schub_{1^n \times w}(x;y) = 
   \sum c_w(a,b,\lambda) \, y_{2+n}^{a_2} \cdots y_{m+n}^{a_m}
   \, (-x_{2+n})^{b_2} \cdots (-x_{m+n})^{b_m} \, s_\lambda(x/y)
\]
where $s_\lambda(x/y)$ is in variables $x_1,\dots,x_{m+n}$ and
$y_1,\dots,y_{m+n}$.

Now restrict to two fixed sets of variables $x_1,\dots,x_N$ and
$y_1,\dots,y_M$, setting $x_i = y_j = 0$ for $i > N$ and $j > M$.
When $n \geq \max(N-1,M-1)$, the only non-zero terms in the above
expression for $\Schub_{1^n \times w}$ are those with all exponents
$a_i$ and $b_i$ equal to zero.  Since this Schubert polynomial is
homogeneous of degree equal to the length of $w$, the partitions
$\lambda$ occurring in these terms all have weight $\ell(w)$.  This
proves:

\begin{thm}
Let $w \in S_{m+1}$ and fix two sets of variables $x_1,\dots,x_N$ and
$y_1,\dots,y_M$.  When $n \geq \max(N-1,M-1)$, the double Schubert
polynomial $\Schub_{1^n \times w}$ in these variables is given by
\[ \Schub_{1^n \times w}(x_1,\dots,x_N,0,\dots,0 ; 
   y_1,\dots,y_M,0,\dots,0) =
   \sum_{\lambda \partof \ell(w)} c_w(0,0,\lambda) \, s_\lambda(x/y)
   \,.
\]
\end{thm}

Comparing with equations (\ref{eqn:stanley}) and (\ref{eqn:stanstab})
we obtain (since $\alpha_{w \lambda} = \alpha_{w^{-1} \lambda'}$):
\begin{cor}
Stanley's coefficient $\alpha_{w \lambda}$ is equal to
$c_w(0,0,\lambda')$.
\end{cor}
Thus the formula (\ref{eqn:schub}) writes a Schubert polynomial as a
symmetric part equal to Stanley's symmetric function plus additional
non-symmetric terms.  For example, if $w = 2\,4\,3\,1$ as in the above
example, we have $F_w(x) = s_{\Pca}(x)$.


The identity $\alpha_{w \lambda'} = \alpha_{w^{-1} \lambda}$ becomes a
special case of the identity $c_{w^{-1}}(b,a,\lambda') =
c_w(a,b,\lambda)$, which in turn follows from the formula
$c_{\mu^\vee}(r^\vee) = c_\mu(r)$ of \cite{buch.fulton:chern}.  Here
$r^\vee$ are the rank conditions obtained by mirroring the rank
diagram for $r = (r_{ij})$ in a vertical line, so $r^\vee_{ij} =
r_{n-j,n-i}$, and $\mu^\vee$ is the sequence $(\mu_n', \dots, \mu_1')$
of conjugate partitions in the opposite order.

\begin{example}
Let $w_0 = m \cdots 2\,1$ be the longest permutation in $S_m$.  Then
we have $r_{w_0}(p,q) = \max(p+q-m, 0)$.  The rectangle diagram
associated to $w_0$ therefore has exactly $i$ non-empty rectangles in
the $i^\Th$ row for $1 \leq i \leq m-1$, and these are centered around
the middle.  All other rectangles are empty.
\[ \raisebox{-1.67cm}{\pic{45}{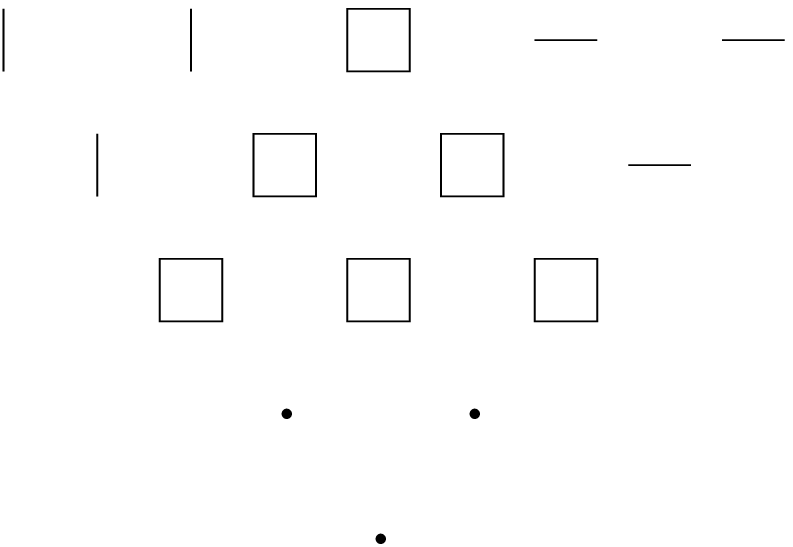}}
   \hspace{-30cm} \left\lbrace \hspace{30cm}
   \begin{matrix} \vspace{1cm} \\ \end{matrix} \right\rbrace m-1
   \smpicskip
\]
\smpicskip
We will use this diagram to compute Stanley's symmetric function for
$w_0$.  The idea is that in order for the algorithm to produce a
contribution to $F_{w_0}$, all boxes must travel north-west until
they meet the last non-empty rectangle in this direction, and from that
point they must travel north-east.

Let $P_r^{(k)}$ denote the element in $\Lambda^{\otimes 2m-3-k}$ given
by this rectangle diagram with the top $k$ rows removed.  In
particular $P_r^{(0)} = P_r$ is the element associated to the whole
diagram.  The terms $s_\lambda$ in Stanley's symmetric function
$F_{w_0}$ are in 1-1 correspondence with terms in $P_r$ of the form
\[ \underbrace{1 \otimes \cdots \otimes 1}_{m-2} \otimes s_\lambda 
   \otimes \underbrace{1 \otimes \cdots \otimes 1}_{m-2} \,.
\]
One may check that, in order for a term
\[ s_{\mu_1} \otimes \cdots \otimes s_{\mu_{2m-3-k}} \]
in $P_r^{(k)}$ to contribute to $F_{w_0}$, the partition $\mu_i$ must
be empty if the $i^\Th$ rectangle in the $k+1^\St$ row of the
rectangle diagram is empty, while it must have length at most one
unless the $i^\Th$ rectangle is the leftmost non-empty rectangle in
the $k+1^\St$ row.  To be precise, the term $s_{\mu_1} \otimes \cdots
\otimes s_{\mu_{2m-3-k}}$ contributes to $F_{w_0}$ only if $\mu_i$ has
length at most one for $i \neq m-1-k$, and is empty when $i \leq
m-2-k$ and when $i \geq m$.  The reason is that all rectangles in the
diagram have height at most one, which means that any $\sigma_i$ in the
algorithm can have length at most one.  So if any $\mu_i$ has two or
more rows, boxes are forced to the right, creating a new partition
with too many rows if $i \neq m-1-k$.

An examination of the algorithm then shows that $P_r^{(k)}$ contains
only one such term, with coefficient 1.  This term has $\mu_{m-1-k}$
equal to the staircase partition with $m-1-k$ rows, $\mu_{m-1-k} =
(m-1-k, \dots, 2, 1)$, while any other non-empty partitions $\mu_i$ is
a single row with $m-1-k$ boxes, $\mu_i = (m-1-k)$.

Taking $k = 0$, we see that $F_{w_0} = s_{(m-1,m-2,\dots,2,1)}$.  This
was first proved by Stanley \cite{stanley:on*1}, and implies that the
number of reduced words for $w_0$ is equal to the number of standard
tableaux on the staircase partition with $m-1$ rows.
\end{example}

\section{Redundant rank conditions and products of permutations}
\lab{sec:redundant}

Suppose we are given a sequence of bundles $E_0 \to E_1 \to \dots \to
E_n$ and a set of rank conditions $r = \{ r_{ij} \}$ for this
sequence.  The degeneracy locus $\Omega_r(E_\bull)$ is then the subset
of points $x \in X$ over which the maps on fibers satisfy all of the
inequalities
\[ \rank(E_i(x) \to E_j(x)) \leq r_{ij} \] 
for $i < j$.  Some of these inequalities may be redundant in the sense
that they follow from other inequalities.  It is easy to see that the
inequality involving the number $r_{ij}$ is redundant if and only
if this number is equal to one of $r_{i,j-1}$ or $r_{i+1,j}$.  In
other words, the inequality involving $r_{ij}$ is necessary if and
only the rectangle $R_{ij}$ is not empty.

Now suppose there are integers $0 \leq p \leq q \leq n$ such that the
rectangle $R_{ij}$ is empty whenever exactly one of $i$ and $j$ is in
the interval $[p,q]$.  In this case the degeneracy locus
$\Omega_r(E_\bull)$ is the (scheme-theoretic) intersection of two
larger loci $\Omega_{r'}(E'_\bull)$ and $\Omega_{r''}(E''_\bull)$ for
the sequences $E'_\bull: E_p \to E_{p+1} \to \dots \to E_q$ and
$E''_\bull: E_1 \to \dots \to E_{p-1} \to E_{q+1} \to \dots \to E_n$,
where $r'$ and $r''$ are the restrictions of the rank conditions $r =
\{r_{ij}\}$ to these sequences.  We will say that $E'_\bull$ is an
independent subsequence.  Note that if $p = q$, the bundle $E_p$ is
redundant and can be removed from the sequence $E_\bull$ without
changing $\Omega_r(E_\bull)$.  This special case was described in
\cite{buch.fulton:chern}.

When $E'_\bull$ is an independent subsequence, the rectangle diagram
for the rank conditions $r'$ simply consists of the rectangles
$R_{ij}$ for $p \leq i < j \leq q$, while the rectangle diagram for
$r''$ contains the remaining non-empty rectangles.
\[ \pic{45}{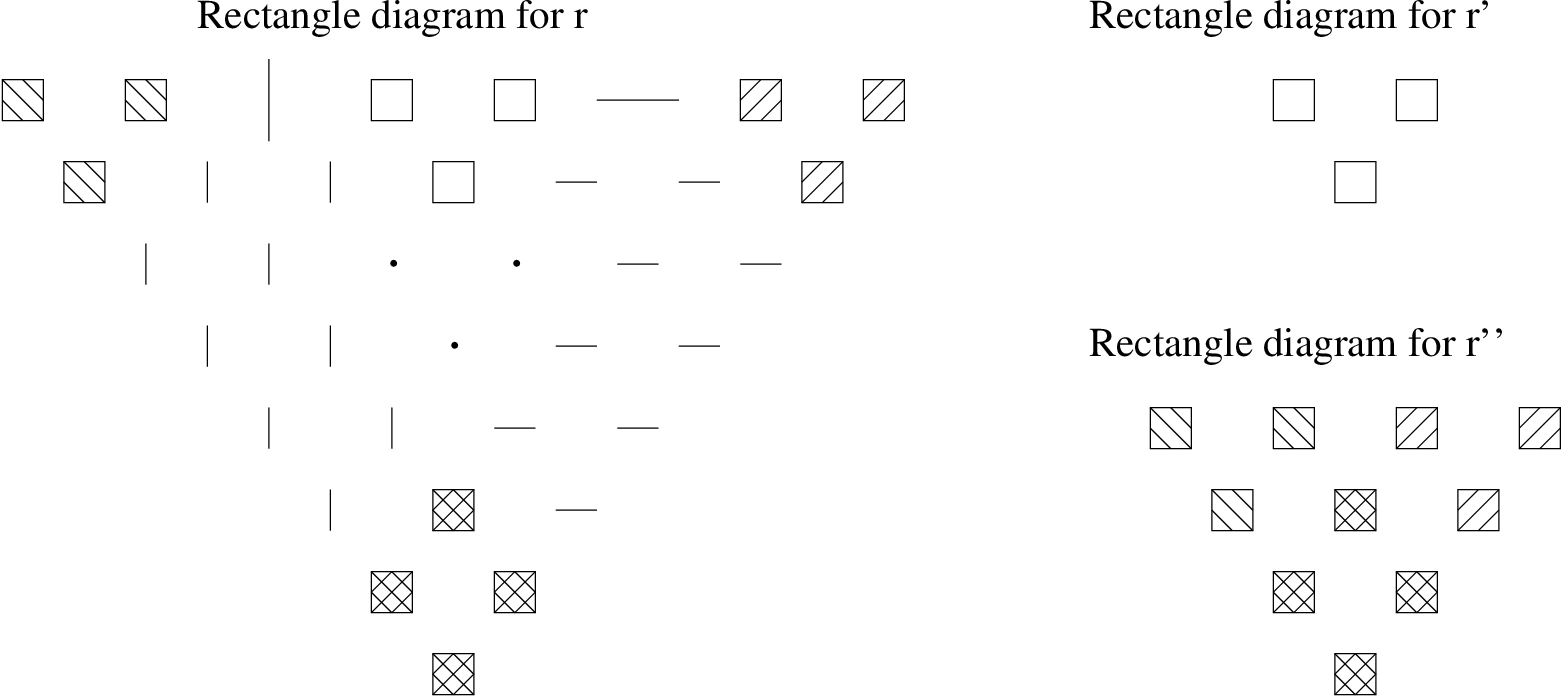} \]
This in particular means that the expected codimensions of
$\Omega_{r'}(E'_\bull)$ and $\Omega_{r''}(E''_\bull)$ add up to that
of $\Omega_r(E_\bull)$.  If all of these loci have their expected
codimensions, then we get the equality $[\Omega_r(E_\bull)] =
[\Omega_{r'}(E'_\bull)] \cdot [\Omega_{r''}(E''_\bull)]$ in the
cohomology ring of $X$.  To see this, note at first that both
$\Omega_{r'}(E'_\bull)$ and $\Omega_{r''}(E''_\bull)$ are
Cohen-Macaulay \cite{lakshmibai.magyar:degeneracy} (see also
\cite[Lemma~A.2]{fulton.pragacz:schubert}).  If $f :
\Omega_{r'}(E'_\bull) \hookrightarrow X$ is the inclusion, we
therefore get
\[ [\Omega_r(E_\bull)]
   = f_* [f^{-1}(\Omega_{r''}(E''_\bull))]
   = f_* f^* [\Omega_{r''}(E''_\bull)]
   = [\Omega_{r'}(E'_\bull)] \cdot [\Omega_{r''}(E''_\bull)] \,.
\]
This means that the formula $P_r$ satisfies
\begin{equation}
P_r = (\underbrace{1 \otimes \dots \otimes 1}_p \otimes P_{r'}
   \otimes \underbrace{1 \otimes \dots \otimes 1}_{n-q}) 
   \cdot \Phi^{q-p+2}_p(P_{r''})
\end{equation}
where multiplication is performed factor-wise, and $\Phi^k_p$ denotes
the $k$-fold coproduct expansion of the $p^\Th$ factor of its
arguments, i.e.\
\begin{multline*}
 \Phi^k_p(s_{\mu_1} \otimes \dots \otimes s_{\mu_p} \otimes \dots
   \otimes s_{\mu_\ell}) = \\
   \sum_{\sigma_1,\dots,\sigma_k} c^{\mu_p}_{\sigma_1,\dots,\sigma_k} \,
   s_{\mu_1} \otimes \dots \otimes s_{\mu_{p-1}} \otimes 
   s_{\sigma_1} \otimes \dots \otimes s_{\sigma_k} \otimes
   s_{\mu_{p+1}} \otimes \dots
   \otimes s_{\mu_\ell} \,.
\end{multline*}

We will apply this to study the Schubert polynomial of a product of
two permutations.  If $w \in S_m$ and $u \in S_n$ are permutations,
define the product $w \times u \in S_{m+n}$ to be the permutation
which maps $i$ to $w(i)$ if $1 \leq i \leq m$, while $m+i$ is mapped
to $m+u(i)$ for $1 \leq i \leq n$.  The rank diagram for this
permutation is equal to that of $1^m \times u$, except the bottom
$2m - 2$ rows are replaced by the rank diagram for $w$.  The diamond
of non-empty rectangles in the rectangle diagram for $w \times u$ is
therefore split into a top part containing the diamond of rectangles
for $u$ and a bottom part with the diamond of rectangles for $w$.
\smpicskip
\[ \pic{70}{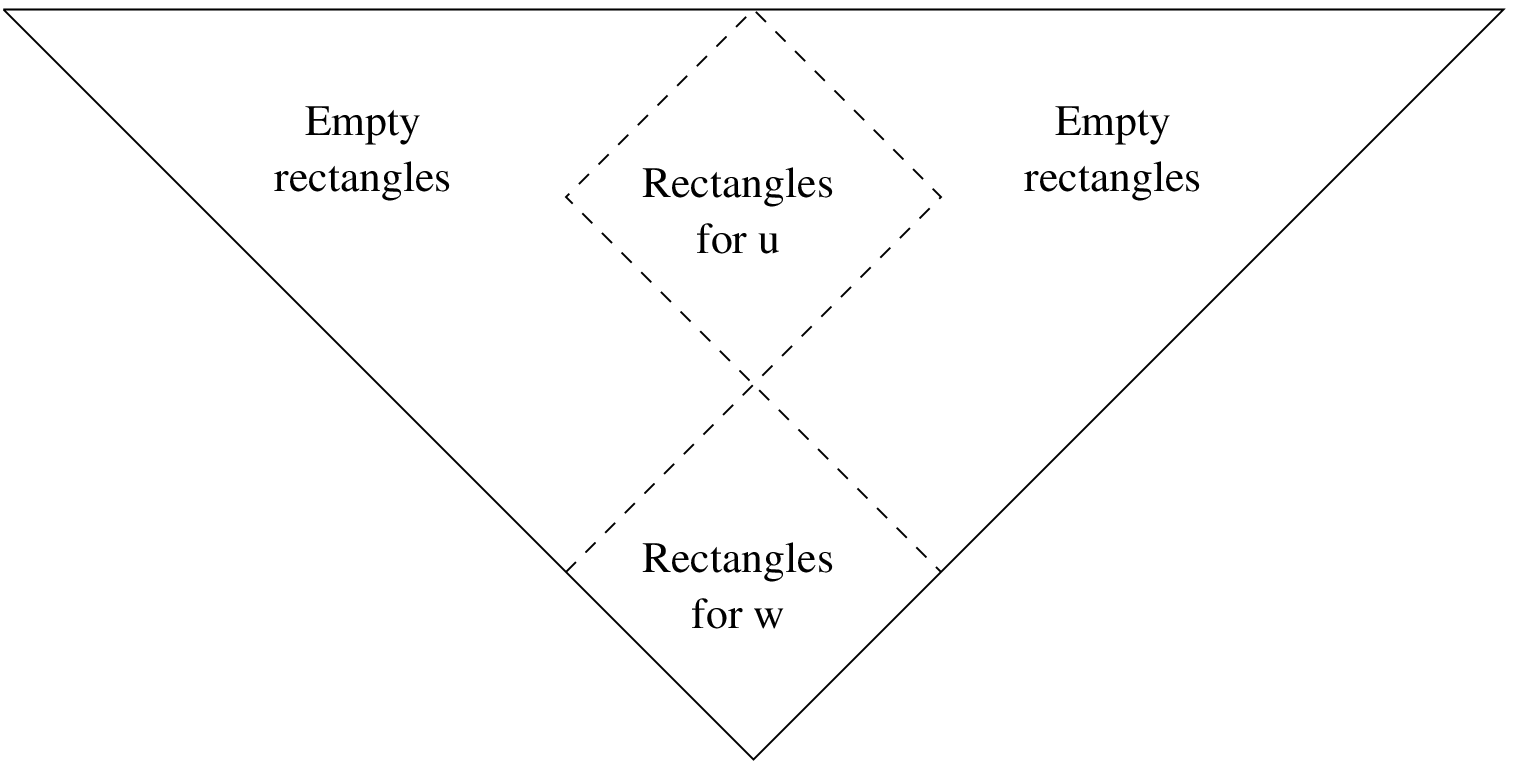} \]
Given a sequence of bundles consisting of a full flag of length
$m+n-1$ followed by a full dual flag of the same length as in Section
\ref{sec:schubert}, we deduce that the locus $\Omega_{w \times u}$ is
the intersection of the loci $\Omega_w$ and $\Omega_{1^m \times u}$.
We therefore recover the well known formula
\cite[(4.6)]{macdonald:notes}
\begin{equation} 
\lab{eqn_prodperm} 
\Schub_{w \times u} = \Schub_w \cdot \Schub_{1^m \times u}
\end{equation}
for the Schubert polynomial of a product of two permutations.  This
immediately implies Stanley's identity $F_{w \times u} = F_w \cdot
F_u$ \cite{stanley:on*1}.  Note that the same argument shows that
\refeqn{eqn_prodperm} also holds for Fulton's universal Schubert
polynomials \cite{fulton:universal}.

\section{Relations to a conjectured Littlewood-Richardson rule}
\lab{sec:conj}

In this final section we will discuss relations with Stanley symmetric
functions of a generalized Littlewood-Richardson rule which is
conjectured in \cite{buch.fulton:chern}.  We will need the notions of
(semistandard) Young tableaux and multiplication of tableaux, see for
example \cite{fulton:young}.

A {\em tableau diagram\/} for a set of rank conditions $r = (r_{ij})$
is a filling of all boxes in the corresponding rectangle diagram with
integers, so that each rectangle $R_{ij}$ becomes a tableau $T_{ij}$.
Furthermore, it is required that the entries of each tableau $T_{ij}$
are strictly larger than the entries in tableaux above
$T_{ij}$ in the diagram, within 45 degree angles.  These are the
tableaux $T_{kl}$ with $i \leq k < l \leq j$ and $(k,l) \neq (i,j)$.

A {\em factor sequence\/} for a tableau diagram with $n$ rows is a
sequence of tableaux $(W_1, \dots, W_n)$, which is obtained as
follows: If $n = 1$ then the only factor sequence is the sequence
$(T_{01})$ containing the only tableau in the diagram.  When $n \geq
2$, a factor sequence is obtained by first constructing a factor
sequence $(U_1,\dots,U_{n-1})$ for the bottom $n-1$ rows of the
tableau diagram, and choosing arbitrary factorizations of the tableaux
in this sequence:
\[ U_i = P_i \cdot Q_i \,. \]
Then the sequence
\[ (W_1,\dots,W_n) = (T_{01} \cdot P_1 \,,\, 
   Q_1 \cdot T_{12} \cdot P_2 \,, \dots ,\,
   Q_{n-1} \cdot T_{n-1,n})
\]
is a factor sequence for the whole tableau diagram.

In \cite{buch.fulton:chern} it is conjectured that the coefficient
$c_\mu(r)$ is equal to the number of different factor sequences $(W_1,
\dots, W_n)$ for any fixed tableau diagram for the rank conditions
$r$, such that $W_i$ has shape $\mu_i$ for each $i$.  This conjecture
has been proved (using an involution of Fomin) when all the rectangles
in the fourth row of the rectangle diagram and below are empty, and no
two non-empty rectangles in the third row are neighbors
\cite{buch:combinatorics}.

If $r = (r_{ij})$ are the rank conditions given by a permutation $w$,
then the conjecture implies that Stanley's coefficient $\alpha_{w
\lambda}$ is equal to the number of different tableaux $W$ of shape
$\lambda'$, for which $(\emptyset, \dots, \emptyset, W, \emptyset,
\dots, \emptyset)$ is a factor sequence.  Thus a proof of the general
conjecture will give a new proof that Stanley's coefficients are
non-negative, as well as an interesting way to compute them.

\begin{example}
Let $w = 2\,1\,4\,3 \dots (2p)\,(2p-1) \in S_{2p}$ for some $p > 0$.
Then the rectangle diagram for $w$ has a $1 \times 1$ rectangle in
the middle of row $4i+1$ for $0 \leq i \leq p-1$.  All other
rectangles are empty.
\[ \pic{30}{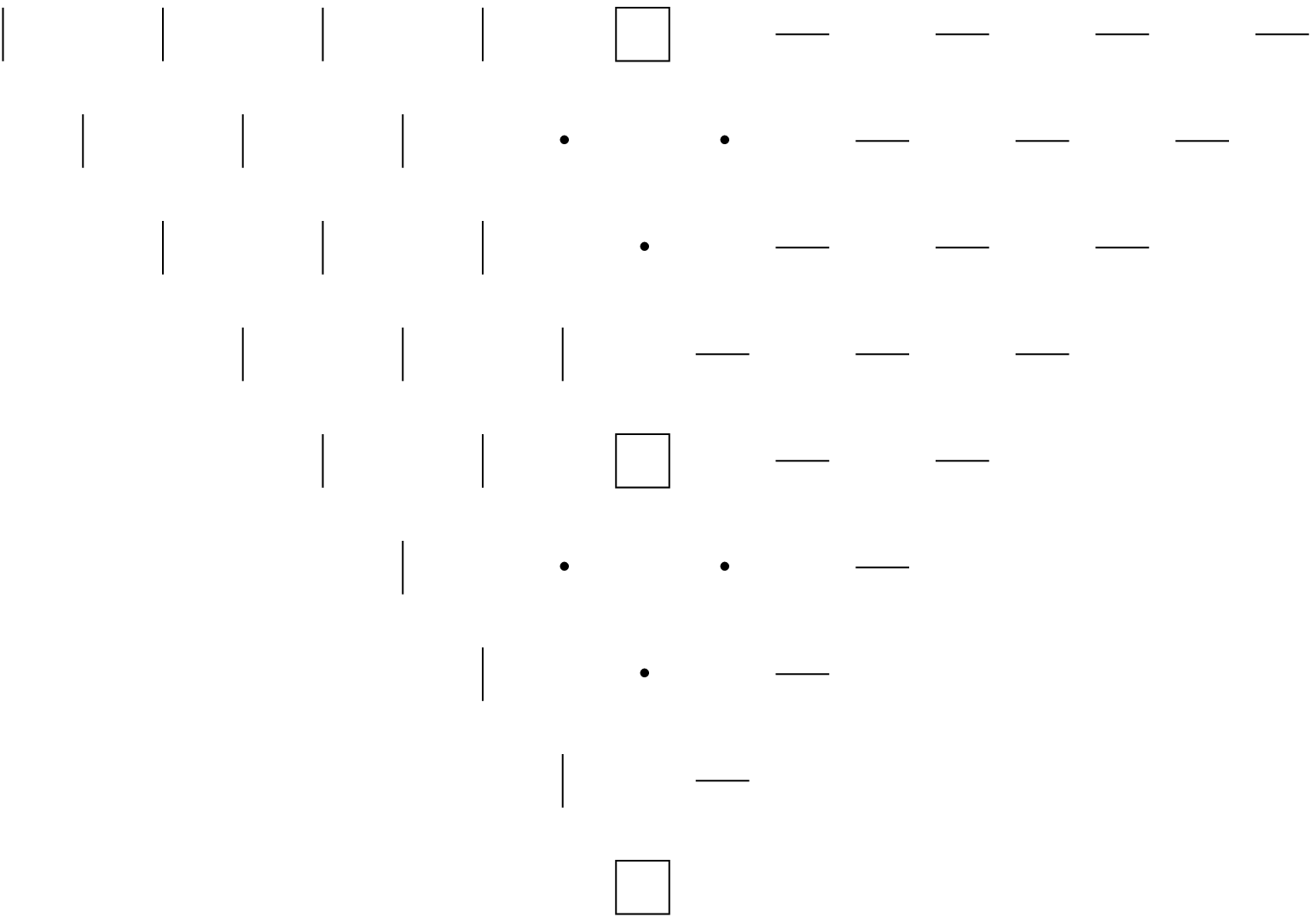} \]
A tableau diagram is obtained by filling the numbers $1, 2, \dots, p$
in these boxes.  It is easy to see that a sequence $(\emptyset, \dots,
\emptyset, W, \emptyset, \dots, \emptyset)$ is a factor sequence for
this diagram if and only if $W$ is a standard tableau with $p$ boxes.
Therefore the conjecture predicts that Stanley's symmetric function is
given by
\[ F_w = \sum_{\lambda \partof p} f^\lambda \, s_\lambda \,. \]
This can be confirmed using Stanley's formula $F_{w \times u} = F_w
\cdot F_u$ \cite{stanley:on*1}.  Let $\sigma = 2\,1 \in S_2$.  Then $w
= \sigma \times \cdots \times \sigma$ (p times), which implies that
\[ F_w = (F_\sigma)^p = (s_{\Pa})^p = \sum_{\lambda \partof p}
   f^\lambda \, s_\lambda \,.
\]
We thank F.~Sottile for showing us a different proof of this fact.
\end{example}

Using a criterion for factor sequences given in
\cite{buch:combinatorics}, one may also prove that the conjectured
Littlewood-Richardson rule gives the correct prediction for Stanley's
symmetric function of a longest permutation $w_0$.

In general, Stanley's symmetric function $F_w$ is known to have a
minimal term $s_{\lambda(w)}$ and a maximal term $s_{\mu(w)}$, both
occurring with coefficient one.  If $w \in S_{m+1}$, define
\[ r_p(w) = \# \{ q \mid q < p \text{ and } w(q) > w(p) \} \]
for $1 \leq p \leq m+1$, and let $\lambda(w)$ be the partition
obtained by arranging the numbers $r_1(w), \dots, r_{m+1}(w)$ in
decreasing order.  Let $\mu(w)$ be the conjugate of the partition
$\lambda(w^{-1})$.  Then $\alpha_{w,\lambda(w)} = \alpha_{w,\mu(w)} =
1$, and any partition $\lambda$ with $\alpha_{w \lambda} \neq 0$ is
between $\lambda(w)$ and $\mu(w)$ in the dominance order
\cite{stanley:on*1}.

Let $\{T_{ij}\}_{1 \leq i < j \leq 2m}$ be a tableau diagram for (the
rank conditions given by) $w$.  There are two extremal ways to form a
factor sequence $(\emptyset, \dots, \emptyset, W, \emptyset, \dots,
\emptyset)$ for this diagram.  The first is to make all factorizations
of inductive factor sequences $(U_1,\dots,U_k)$ be ``rightward''
whenever possible.  This means that when factoring $U_i$ into $U_i =
P_i \cdot Q_i$, we take $P_i = \emptyset$ and $Q_i = U_i$ for $i \neq
m$ while we take $P_m = U_m$ and $Q_m = \emptyset$ (if $k \geq m$).
The middle tableau in the final factor sequence then is
\[ W_{\text{right}} = T_m \cdot T_{m+1} \cdot \ldots \cdot T_{2m-1} \]
where
\[ T_j = T_{0j} \cdot T_{1j} \cdot \ldots \cdot T_{m-1,j} \,. \]
Note that each tableau $T_j$ has only one column.  If we set $p =
2m+1-j$ and $q = w^{-1}(i+2)$ then $T_{ij}$ is non-empty if and only
if $q < p$ and $w(q) > w(p)$.  It follows that $T_j$ has exactly
$r_p(w)$ boxes.

We claim that $W_{\text{right}}$ has shape $\lambda(w)'$,
corresponding to the maximal term of $F_{w^{-1}}$.  It is enough to
show that the if $T_l$ and $T_j$ both have a box in row $t$ and $l <
j$, then the box in $T_l$ is smaller than the one in $T_j$.  To prove
this, let the $t^\Th$ box in $T_l$ come from $T_{kl}$ and the $t^\Th$
box in $T_j$ come from $T_{ij}$.  If the box in $T_{kl}$ is not
smaller than the box in $T_{ij}$ then $k < i$.  Now since the tableau
$T_{il}$ must be as wide as $T_{kl}$ and as tall as $T_{ij}$, this
tableau $T_{il}$ can't be empty.  Similarly, if $T_{hj}$ corresponds
to a box over $T_{ij}$ in $T_j$, then $T_{hl}$ gives a corresponding
box in $T_l$.  This shows that the boxes corresponding to $T_{kl}$ and
$T_{ij}$ in $T_l$ and $T_j$ was not in the same row, a contradiction.

Similarly one can show that the tableau obtained by ``leftward''
factorizations,
\[ W_{\text{left}} = 
   (T_{0,m} \cdot T_{0,m+1} \cdots T_{0,2m-1}) \cdot
   (T_{1,m} \cdots T_{1,2m-1}) \cdots
   (T_{m-1,m} \cdots T_{m-1,2m-1}) \,,
\]
has shape $\mu(w)'$.



\providecommand{\bysame}{\leavevmode\hbox to3em{\hrulefill}\thinspace}


\end{document}